\documentclass{article}
\usepackage{amsmath,amstext,amsthm,amsfonts}
\usepackage[dvips]{graphicx}
\usepackage{epic,eepic}

\theoremstyle{plain}
\newtheorem{theorem}{Theorem }[section]
\newtheorem{proposition}[theorem]{Proposition}
\newtheorem{lemma}[theorem]{Lemma}
\newtheorem{corollary}[theorem]{Corollary}
\newtheorem{maintheorem}{Theorem}
\newtheorem{theorem A}{Theorem A}
\newtheorem{theorem B}{Theorem B}
\newtheorem{theorem C}{Theorem C}
\newtheorem{theorem D}{Theorem D}

\theoremstyle{definition}
\newtheorem{remark}[theorem]{Remark}

\newtheorem{definition}[theorem]{Definition}

\newtheorem{question}{Question}
\newtheorem{conjecture}{Conjecture}
\newcommand{\field}[1]{\mathbb{#1}}
\newcommand{\real}{\field{R}}

\renewcommand{\natural}{\field{N}}
\newcommand{\rational}{\field{Q}}

\newcommand{\al} {\alpha}       
        
\newcommand{\ga} {\gamma}    \newcommand{\Ga}{\Gamma}
       
\newcommand{\ep} {\epsilon}
\newcommand{\vep}{\varepsilon}

\newcommand{\si} {\sigma}

\newcommand{\vfi}{\varphi}

\newcommand{\SC}{{\cal C}}

\newcommand{\SI}{{\cal I}}

\newcommand{\SL}{{\cal L}}
\newcommand{\SM}{{\cal M}}

\newcommand{\SR}{{\cal R}}

\begin{document}

\title{Geometrical \emph{versus} Topological Properties of Manifolds }

\author{Carlos Matheus\footnote{C. Matheus is supported by Faperj/Brazil}
\text{ and}
Krerley Oliveira\footnote{K. Oliveira is supported by Pronex-CNPq/Brazil and Fapeal/Brazil} }

\date{October 30, 2003}

\maketitle
\begin{abstract} Given a compact $n$-dimensional immersed Riemannian manifold $M^n$ in some Euclidean space  we prove that if
the Hausdorff dimension of the singular set of the Gauss map is small, then  $M^n$ is homeomorphic to the sphere
$S^n$.

Also, we define a concept of finite geometrical type and prove that finite
geometrical type hypersurfaces with small set of points of zero Gauss-Kronecker
curvature are topologically the sphere minus a finite number of
points. A characterization of the $2n$-catenoid is obtained.
\end{abstract}

\section{Introduction}

   Let $f:M^n\rightarrow N^m$ be a $C^1$ map. We denote by

$$rank(f):=\min\limits_{p\in M}rank(D_p f).$$

 If $n=\dim M=\dim N=m$, let $C:=\left\{ p\in M: \det D_p f = 0 \right\}$ the set of \emph{critical points}
of $f$ and $S:=f(C)$ the set of \emph{critical values} of $f$.

Now, let $M^n$ a compact, connected, boundary less,
$n$-dimensional manifold. Denote by $H_s$ the $s$-dimensional
Hausdorff measure and $\dim_H(A)$ the Hausdorff dimension of
$A\subset M^n$. For definitions see section~\ref{notacoes} below.
Let $x$ be an immersion $x:M^n\rightarrow \real^{n+1}$. In this case,
let $G:M^n\rightarrow S^n$ the Gauss map associated to $x$, $C$
the critical points of $G$ and $S$ the critical values of $G$. We
denote by $dim_H(x):=dim_H(S)$.  By Moreira's improvement of
Morse-Sard theorem (see~\cite{Moreira}), since $G$ is a smooth
map, we have that $\dim_{H}(S)\leq n-1$.

In other words, if $\SI mm=\{ x:M\rightarrow
\real^{n+1} : x \text{ is an immersion} \}$, then $\sup\limits_{x\in\SI mm} dim_{H}(x)\leq n-1$. Clearly, this
supremum could be equal to $n-1$, as some immersions of $S^n$ in $\real^{n+1}$
show (e.g., immersions with ``cylindrical pieces''). Our interest here is the number $\inf dim_H(x)$.
Before discuss this, we introduce some definitions.

\begin{definition}\label{rank de imersao}Given an immersion
$x:M^n\rightarrow\real^{n+1}$ we define $rank(x):=rank(G)$, where
$G$ is the Gauss map for $x$.
\end{definition}

\begin{definition}   We denote by $\SR(k)$ the set
$\SR(k)=\{x\in \SI mm : rank(x)\geq k\}.$ Define by $\alpha_k(M)$ the numbers:

$$\al_k(M)=\inf\limits_{x\in\SR(k)}dim_H (x), k=0,\dots,n$$

If $\SR(k)=\emptyset$ we define $\alpha_k(M) =n-1$.

\end{definition}

     Now, we are in position to state our first result:

\begin{maintheorem}\label{theoremA}If $M^n$ is a compact manifold with $n\geq 3$ such that
$\al_k(M^n)<k-[\frac{n}{2}]$, for some integer $k$, then $M^n$ is homeomorphic to $S^n$ ($[r]$ is the integer part of $r$).
\end{maintheorem}

 The proof of this theorem in the cases $n=3$ and $n\geq 4$ are quite different. For higher dimensions, we can use the generalized Poincar\'e
Conjecture (Smale and Freedman) to obtain that the given manifold is a sphere.
Since the Poincar\'e Conjecture is not available in three dimensions, the proof,
in this case, is a little bit different. We use a characterization theorem due
to Bing to compensate the loss of Poincar\'e Conjecture, as commented before.

   To prove this theorem in the case $n=3$, we proceed as follows:

   \begin{itemize}
   \item  By a theorem of Bing (see \cite{Bing}), we just need to prove that every piecewise smooth simple curve $\gamma$
   in $M^3$ lies in a topological cube $\SR$ of $M^3$;

   \item In order to prove it, we shall show that it is enough to prove for $\gamma \subset
   M-G^{-1}(S)$ and that $G:M-G^{-1}(S)\rightarrow S^3-S$ is a diffeomorphism;

   \item Finally, we produce a cube $\tilde{\SR}\supset G(\gamma)$ in $S^3-S$ and we
   obtain $\SR$ pulling back this cube by $G$

   \end{itemize}
Observe that by~\cite{Cohen}, in three dimensions always there are Euclidean codimension 1 immersions.
      In particular, it is reasonable to consider the following consequence of the Theorem A:

\begin{corollary}\label{conjecture1} The following statement is
equivalent to Poincar\'e Conjecture : ``Simply connected $3$-manifolds
admits Euclidean codimension one immersions with rank at least 2 and Hausdorff dimension of the singular set for
this Gauss map less than 1''.
\end{corollary}

   For a motivation of this conjecture and some comments about three dimensional manifolds see the section
   \ref{finalremarks}.

Our motivation behind proving this theorem are results by do Carmo,
Elbert~\cite{dCE} and Barbosa, Fukuoka, Mercuri~\cite{BFM}.
Roughly speaking, they obtain topological results about certain
manifolds provided that there are special codimension 1 immersions
of them. These results motivate the question : how the space of
immersions (extrinsic information) influences the topology of $M$
(intrinsic information)? The Theorems A and B below are a partial answer
to this question. The proofs of the theorems depends on the concept of Hausdorff
dimension. Essentially, Hausdorff dimension is a fractal dimension
that measures how ``small'' is a given set with respect to usual
``regular'' sets (e.g., smooth submanifolds, that always have
integer Hausdorff dimension).

In section 6 of this paper we obtain the following generalizations of Theorem~\ref{theoremA} and~\ref{theoremB}:

\begin{definition}Let $\overline{M}^n$ a compact (oriented) manifold and
$p_1,\dots,p_k\in\overline{M}^n$. Let $M=\overline{M}^n-\{p_1,\dots,p_k\}$. An
immersion $x:M^n\rightarrow \real^{n+1}$ is of \emph{finite geometrical type} (in a weaker sense than that of~\cite{BFM}) if $M^n$
is complete in the induced metric, the Gauss map $G:M^n\rightarrow S^n$ extends
continuously to a function $\overline{G}:\overline{M}^n\rightarrow S^n$ and the set
$G^{-1}(S)$ has $H_{n-1}(G^{-1}(S))=0$ (this last condition occurs
if $rank(x)\geq k$ and $H_{k-1}(S)=0$).
\end{definition}

The conditions in the previous definition are satisfied by complete
hypersurfaces with finite total curvature whose Gauss-Kronecker curvature $H_n=
k_1\dots k_n$ does not change of sign and vanish in a small set, as showed by \cite{dCE}. Recall that
a hypersurface $x:M^n\rightarrow \real^{n+1}$ has total finite curvature if
$\int_M |A|^n dM<\infty$, $|A|=(\sum\limits_i k_i^2)^{1/2}$, $k_i$ are the
principal curvatures. With these observations, one has :

\begin{maintheorem}\label{theoremB}If $x:M^n\rightarrow\real^{n+1}$ is a hypersurface with finite
geometrical type and
$H_{k-[\frac{n}{2}]}(S)=0$, $rank(x)\geq k$. Then $M^n$ is topologically a sphere minus a finite number
of points, i.e., $\overline{M}^n\simeq S^n$. In particular, this result holds for
complete hypersurfaces with finite total curvature and
$H_{k-[\frac{n}{2}]}(S)=0$, $rank(x)\geq k$.
\end{maintheorem}

For even dimensions, we follow \cite{BFM} and improve Theorem~\ref{theoremB}.
In particular,
we obtain the following  characterization of $2n$-catenoids, as the unique minimal hypersurfaces of finite
geometrical type.

\begin{maintheorem}\label{theoremC} Let $x:M^{2n}\rightarrow\real^{2n+1}$, $n\geq 2$ an immersion of finite
geometrical type with $H_{k-n}(S)=0$, $rank(x)\geq k$. Then $M^{2n}$ is
topologically a sphere minus two points. If $M^{2n}$ is minimal, $M^{2n}$ is a
$2n$-catenoid.
\end{maintheorem}

\section{Notations and Statements}\label{notacoes}

Let $M^n$ be a smooth manifold. Before starting the proofs of the statements we fix some
notations and collect some (useful) standard propositions about Hausdorff
dimension (and limit capacity, another fractal dimension). For the proofs of
these propositions we
refer~\cite{Falconer}.

Let $X$ a compact metric space and $A\subset X$. We define the $s$-dimensional
Hausdorff measure of $A$ by

$$H_s(A):=\lim\limits_{\vep\rightarrow 0} \inf\{
\sum\limits_{i}(diam\ U_i)^s : A\subset\bigcup U_i, U_i \text{ is open }\newline
\text{and }diam (U_i)\leq\vep, \forall i\in\natural \}.$$
The \emph{Hausdorff
dimension} of $A$ is $\dim_H(A):=\sup\{d\geq 0:H_d(A)=\infty\}=\inf\{d\geq 0:
H_d(A)=0\}.$ A remarkable fact is that $H_n$ coincides with Lebesgue
measure for a smooth manifold $M^n$.

 A related notion are the lower and upper \emph{limit capacity} (sometimes called box
counting dimension) defined by
$$\underline{\dim_B}(A):=\liminf\limits_{\vep\rightarrow 0}\log n(A,\vep)/ (-
\log\vep),\ \overline{\dim_B}(A):=\limsup\limits_{\vep\rightarrow 0}\log
n(A,\vep)/(-\log\vep),$$ where $n(A,\vep)$ is the minimum number of $\vep$-balls
that cover $A$. If $d(A)=\underline{\dim_B}(A)=\overline{\dim_B}(A)$, we say that the
limit capacity of $A$ is $\dim_B(A)=d(A)$.

These fractal dimensions satisfy the properties expected for ``natural'' notions
of dimensions. For instance, $\dim_H(A)=m$ if $A$ is a smooth
$m$-submanifold.

\begin{proposition}\label{fatos basicos}The properties listed below hold :

\begin{enumerate}

\item $\dim_H(E)\leq\dim_H(F)$ if $E\subset F$;

\item  $\dim_H(E\cup F)=\max\{\dim_H(E),
\dim_H(F)\}$;

\item  If $f$ is a Lipschitz map with Lipschitz constant $C$, then $H_s(f(E))\leq C\cdot H_s(E)$.
As a consequence,  $\dim_H(f(E))\leq\dim_H E$;

\item If $f$ is a bi-Lipschitz map (e.g.,
diffeomorphisms), $\dim_H(f(E))=\dim_H(E)$;

\item $dim_H(A) \leq \underline{dim_B}(A)$.

\end{enumerate}

Analogous properties holds for lower
and upper limit capacity. If $E$ is countable, $\dim_H(E)=0$ (although we may
have $\dim_B(E)>0$).

\end{proposition}

When we are dealing with product spaces, the relationship between Hausdorff
dimension and limit capacity are the product formulae :

\begin{proposition}\label{produto}
 $\dim_H(E)+\dim_H(F)\leq\dim_H(E\times F)\leq\dim_H(E)+
\overline{\dim_B}(F)$. Moreover, $c\cdot H_s(E)\cdot H_t(F)\leq H_{s+t}(E\times F)\leq
C\cdot H_s(E)$,
where $c$ depends only on $s$ and $t$, \ $C$ depends only on $s$ and $\overline{
\dim_B}(F)$.
\end{proposition}

Before starting the necessary lemmas to prove the central results, we observe that it follows from lemma above that if $M$ and $N$ are diffeomorphic $n$-manifolds then $\al_k(M)=\al_k(N)$. This
proves :

\begin{lemma}\label{constante alfa}  The numbers
$$\al_k(M)=\inf\limits_{x\in\SR(k)}dim_H (x), \text{ for }k=0,\dots,n$$
are smooth invariants of $M$.
\end{lemma}

In particular, if $n=3$ we also have that $\alpha_k$ are topological invariants. It is a consequence of a theorem due to Moise~\cite{Moise}, which state that if $M$ and
$N$ are homeomorphic $3$-manifolds then they are diffeomorphic. Then, the following conjecture arises from the Theorem~ \ref{theoremA}

\begin{conjecture}\label{principio variacional}\textit{If $M^3$ is simply
connected, then  $$\alpha_2(M^3)=\inf\limits_{x\in\SR(2)} dim_{H}(x)< 1$$}
\end{conjecture}

R. Cohen's theorem~(\cite{Cohen}) says that there are immersions of compact
$n$-manifolds $M^n$ in $\real^{2n-\al(n)}$ where $\al(n)$ is the number of 1's in the binary expansion
of $n$. This implies, for the case $n=3$, that we always have
that $\SI mm\neq\emptyset$. In particular, the implicit hypothesis of existence
of codimension 1 immersions in Theorem~\ref{theoremA} is not too restrictive and our conjecture is
reasonable. We point out
that conjecture~\ref{principio variacional} is true if Poincar\'e conjecture
holds and, in this case, $\sup\limits_{x\in\SI mm}rank(x)=3$ and
$\inf\limits_{x\in\SR(k)} dim_{H}(x)=0, \text{ for all } 0\leq k\leq 3$. A corollary of
the theorem~\ref{theoremA} and this observation is:

\begin{corollary}\label{Poincare} The Poincar\'e Conjecture is equivalent to the conjecture~\ref{principio variacional}.
\end{corollary}

From this, a natural approach to conjecture~\ref{principio variacional} is a deformation and desingularization
argument for metrics given by pull-back of immersions in $\SI mm$. We observe
that Moreira's theorem give us $\al_2(M^3)\leq 2$. This motivates the following
question, which is a kind of step toward Poincar\'e Conjecture. However, this
question is of independent interest, since it can be true even if Poincar\'e
Conjecture is false :
\begin{question}\textit{For simply connected $3$-manifolds, is true that
$\al_2(M^3)<2$ ? }
\end{question}

\section{Some lemmas}

In this section, we prove some useful facts on the way to establish the Theorems~\ref{theoremA}, B. The first one relates the Hausdorff dimension of
subsets of smooth manifolds and rank of smooth maps :

\begin{proposition}\label{dim pre imagem} Let $f:M^m\rightarrow N^n$ a $C^1$-map and $A\subset N$.
Then $\dim_H f^{-1}(A)$ $\leq \dim_H (A) + n - rank(f)$.
\end{proposition}

\begin{proof} The computation of Hausdorff dimension is a local problem. So, we can
consider $p\in f^{-1}(S)$, coordinate neighborhoods $p\in U$,
$f(p)\in V$ fixed and $f=(f_1,\dots,f_n):U\rightarrow V$. Making a change of
coordinates (which does not change Hausdorff dimensions), we can suppose that $\widetilde{f}=(f_1,\dots,f_r)$ is a
submersion, where $r=rank(f)$. By the local form of submersions, there is
$\vfi$ a diffeomorphism such that $\widetilde{f}\circ\vfi(y_1,\dots,y_m)=(y_1,\dots,
y_r)$. This implies that $f\circ\vfi(y_1,\dots,y_m)=(y_1,\dots,y_r,
g(\vfi(y_1,\dots,y_m))$. Then, if $\pi$ denotes the projection in the $r$ first
variables, $x\in f^{-1}(S)\Rightarrow \pi\vfi^{-1}(x)\in\pi(S)$, i.e.,
$f^{-1}(S)\subset\vfi(\pi(S)\times\real^{n-r})$. By properties of Hausdorff
dimension (see section 2), we have $dim_H f^{-1}(S)\leq dim_H(\pi(S)\times
\real^{n-r})\leq dim_H\pi(S)+ \overline{dim_B}(\real^{n-r})\leq dim_H(S)+n-r$.
This concludes the proof.
\end{proof}

The second proposition relates Hausdorff dimension with topological results.

\begin{proposition}\label{transv}Let $n\geq 3$ and $F$ is a closed
subset of a $n$-dimensional connected (not necessarily compact) manifold $M^n$.
If the Hausdorff dimension of $F$ is strictly less than $n-1$ then $M^n-F$ is
connected. If $M^n=\real^n$ or $M^n=S^n$, $F$ is compact and the Hausdorff
dimension of $F$ is strictly less than $n-k-1$ then $M^n-F$ is $k$-connected
(i.e., its homotopy groups $\pi_i$ vanishes for $i\leq k$).
\end{proposition}

\begin{proof}First, if $F$ is a closed subset of $M^n$ with
Hausdorff dimension strictly less than $n-1$, $x,y\in M^n-F$, take
$\ga$ a path from $x$ to $y$ in $M^n$. Since $n\geq 3$, we can
suppose $\ga$ a smooth simple curve (by transversality). In this
case, $\ga$ admits some compact tubular neighborhood $\SL$. For
each $p\in\ga$, denote $\SL_p$ the $\SL$-fiber passing through
$p$. By hypothesis, $dim_{H}(F\cap \SL_p)<n-1$ $\forall p$. In
this case, the tubular neighborhood $\SL$ is diffeomorphic to
$\ga\times D^{n-1}$, the fibers $\SL_p$ are $p\times D^{n-1}$
($D^{n-1}$ is the $(n-1)$-dimensional unit disk centered at $0$)
and $\ga$ is $\ga\times{0}$. Then, since $F$ is closed, it is easy
that every $x\in\ga$ admits a neighborhood $V(x)$ such that for
some sequence $v_n=v_n(x)\rightarrow 0$ holds
$(V(x)\times{v_n})\cap F=\emptyset$. Moreover, again by the fact
that $F$ is closed, any vector $v$ sufficiently close to some
$v_n$ satisfies $(V(x)\times{v})\cap F=\emptyset$. With this in
mind, by compactness of $\ga$, we get some finite cover of $\ga$
by neighborhoods as described before. This guarantees the
existence of $v_0$ arbitrarily small such that
$(\ga\times{v_0})\cap F= \emptyset$. This implies that $M-F$ is
connected.

Second, if $F$ is a compact subset of $M^n=\real^n$, $dim_H F<n-k-1$, let
$[\Ga]\in\pi_i(\real^n-F)$
a homotopy class for $i\leq k$. Choose a smooth representative $\Ga\in[\Ga]$. Define $f:\Ga\times F\rightarrow S^{n-1}$, $f(x,y):=
(y-x)/||y-x||$. We will consider in $\Ga\times F$ the sum norm, i.e., if
$p,q\in\Ga\times F$, $p=(x,y), q=(z,w)$ then $||p-q||:=||x-z||+||y-w||$. For
this choice of norm we have
$$||f(p)-f(q)||=\frac{1}{||y-x||\cdot||z-w||}\cdot \Big\|\big\{ (y-x)\cdot||z-w|| +
||y-x||\cdot(z-w) \big\}\Big\| \Rightarrow
$$
$ ||f(p)-f(q)||\leq\frac{\big\| (y-x)\cdot||z-w|| -
||z-w||\cdot(w-z)\big\|}{||y-x||\cdot||z-w||} \ + \ \frac{
\big\|||z-w||\cdot(w-z) -
||y-x||\cdot(w-z)\big\|}{||y-x||\cdot||z-w||}\cdot\Rightarrow
$
$$||f(p)-f(q)||\leq \frac{1}{||y-x||}\cdot\Big\{||(z-x)+(y-w)||\Big\} +
\frac{1}{||y-x||}\cdot\big\arrowvert\{||(z-w)||-||(y-x)||\}\big\arrowvert\Rightarrow$$
$$||f(p)-f(q)||\leq 2\cdot C\cdot ||p-q||
$$
where $C=1/d(\Ga,F)$. We have $d(\Ga,F)>0$ since
these are compact disjoint sets. This computation shows that $f$ is Lipschitz.

Then, we have (Proposition~\ref{fatos basicos},~\ref{produto}) $dim_H f(\Ga\times F)\leq dim_H(\Ga\times F)\leq
\overline{dim_B}(\Ga)+dim_H(F)<i+n-k-1\leq n-1\Rightarrow \exists\ v\notin f(\Ga\times
F)$. Now, $F$ is compact implies that there is a real $N$ such that $F\subset
B_N(0)$. Then, making a translation of $\Ga$ at $v$ direction, we can put, using
this translation as homotopy, $\Ga$
outside $B_N$. Since $\real^n-B_N$ is $n$-connected (for $n\geq 3$), $\pi_i(\real^n-F)=0$.
This concludes the proof.
\end{proof}

\begin{remark}\label{contra-ex. a transv}We remark that the
hypothesis $F$ is closed in the previous proposition is necessary. For example,
take $F=\rational^n$, $M^n=\real^n$. We have $dim_H(F)=0$ ($F$ is a countable
set) but $M^n-F$ is not connected.
\end{remark}

We can think of Proposition~\ref{transv} as a weak type of transversality. In fact,
if $F$ is a compact ($n-2$)-submanifold of $M^n$ then $M-F$ is connected and if
$F$ is a compact ($n-3$)-submanifold of $\real^n$ (or $S^n$) then $\real^n-F$ is
simply connected. This follows from basic transversality. However, our previous
proposition does not assume regularity of $F$, but allows us to conclude the same
results. It is natural these results are true because Hausdorff dimension
translates the fact that $F$ is, in some sense, ``smaller'' than a
$(n-1)$-submanifold $N$ which has optimal dimension in
order to disconnect $M^n$.

For later use, we generalize the first part of
Proposition~\ref{transv} as follows :

\begin{lemma}\label{difeo}Suppose that $\Ga\in\pi_i(M^n)$ is Lipschitz (e.g.,
if $i=1$ and $\Ga$ is a piecewise smooth curve) and let
$K\subset M^n$ compact, $dim_H K<n-i$. Then there are diffeomorphisms $h$ of $M$,
arbitrarily close to identity map, such that $h(\Ga)\cap K=\emptyset$. In particular,
if $[\Ga]\in\pi_i(M^n)$ a homotopy class, $K\subset M^n$ a compact set,
$dim_H(K)<n-i$, there is
a smooth representative $\Ga\in [\Ga]$ such that
$\Ga\cap K=\emptyset$, i.e., $\Ga\in\pi_i(M^n -K)$.
\end{lemma}

\begin{proof}First, consider
a parametrized neighborhood $\phi:U\rightarrow B_3(0) \subset
\real^n$ and  suppose that $\Gamma$ lies in $\overline{V_1}$, where $V_1=\phi^{-1}(B_1(0))$. Let
$K_1 = \phi(K) \subset \real^n$ and $\Gamma_1 = \phi(\Gamma) \subset \real^n$. Consider the map:
$$ F: \Gamma_1\times K_1  \rightarrow \real^n \text{ defined by } F(x,y)= x-y$$
Observe that, since $\Gamma$ is Lipschitz and $\phi$ is a diffeomorphism,
$\overline{dim_B} \Gamma=
\overline{dim_B} \Gamma_1 \leq i$. This implies that $dim_H(F(\Gamma_1 \times K_1)) < n$, since
$dim_H(K)<n-i$. This implies, in particular, that  $\real^n - F(\Gamma_1 \times K_1)$ is an open and dense
subset, since $K$ is compact. Then, we may choose a vector $v \in \real^n - F(\Gamma_1 \times K_1)$ arbitrarily close to 0 such
$(\Gamma_1 + v) \subset B_2(0)$. Since, $ v \in \real^n - F(\Gamma_1 \times K_1)$ we have that
$(\Gamma_1 + v) \cap K_1 = \emptyset$.

To construct $h$ we consider a bump function $\beta:\real^n\rightarrow
[0,1]$, such that $\beta(x)=1$
  if $x\in B_1(0)$ and $\beta(x)=0$  for every $x\in \real^n-B_2(0).$ It is easy to see that
  $h$ defined by:

  $$h(y)=y \text{ if } x\in M-U \text{ and } h(y)=\phi^{-1}(\beta(\phi(y))v + \phi(y)),$$ is a
  diffeomorphism that satisfies $h(\Ga)\cap K =\emptyset$, since $(\Ga_1 + v) \cap K_1 =
  \emptyset$.

  In the general case, we proceed as
follows : first, considering a finite number of parametrized neighborhoods $\phi_i: U_i
\rightarrow B_3(0), i\in\{1,\dots,n\}$  and $V_i = \phi_i^{-1}(B_1(0))$ covering
$\Gamma$, by the previous case,
there exists $h_1$ arbitrarily close to the identity such $h_1(\Gamma) \subset
\bigcup\limits_{i=1}^n V_i $ and such that $h_1(\Gamma\cap \overline{V_1})\cap K
=\emptyset$. Observe that,
 $d(h_1(\Gamma\cap \overline{V_1}), K)>\epsilon_1>0$, since $h_1(\Gamma\cap \overline{V_1})$ is
 a compact set.

The next step is to repeat the previous argument considering $h_2$ arbitrarily close to the
identity, in such way that $h_2(h_1(\Gamma)\cap V_2)\cap K=\emptyset$ and
$h_2(h_1(\Ga))\subset
   \bigcup\limits_{i=1}^n V_i$. If $d(h_2,id)<\frac{\epsilon_1}{2}$ then
   $h_2(h_1(\Ga)\cap
   V_1)\cap K = \emptyset$. Repeating this argument by induction, we obtain that
   $h=h_n\circ
   \dots \circ h_1$ is a diffeomorphism such that $h(\Ga)\cap K =\emptyset$. This
   concludes the
   proof.
\end{proof}

\section{Proof of Theorem A in the case $n=3$}
Before giving a proof for theorem~\ref{theoremA}, we mention a lemma
due to Bing~\cite{Bing} :

\begin{lemma}[Bing]\label{Bing}A compact, connected, $3$-manifold $M$ is topologically
$S^3$ if and only if each piecewise smooth simple closed curve in $M$ lies in a topological cube in $M$.
\end{lemma}

A modern proof of this lemma can be found in~\cite{Rolfsen}. In modern language,
Bing's proof shows that the hypothesis above imply that
\emph{Heegaard} splitting of $M$ is in two balls. This is sufficient to conclude the result.

In fact, Bing's theorem is not stated in~\cite{Bing},~\cite{Rolfsen} as above.
But the lemma holds. Actually, to prove that $M$ is homeomorphic to $S^3$,
Bing uses only that, if a triangulation of $M$ is
fixed, every simple \emph{polyhedral} closed curve lies in a topological cube.
Observe that \emph{polyhedral} curves are piecewise smooth curves, if we choose
a smooth triangulation (smooth manifolds always can be smooth triangulated,
see~\cite{Thurston}, page $194$; see also \cite{Whitney}, page 124).

\begin{proof}[Proof of Theorem A in the case $n=3$]
If $\al_2(M)<1$, there is an immersion $x:M^3\rightarrow \real^4$ such that
$rank(x)\geq 2, dim_H(x)<1$. Let G the Gauss map associated to $x$. By
Propositions~\ref{transv},~\ref{dim pre imagem}, since $dim_H(S)<1$,
$M-G^{-1}(S), S^3-S$ are connected manifolds. Consider $G:M-G^{-1}(S)\rightarrow S^3-S$. This is a proper map
between connected manifolds whose Jacobian never
vanishes. So it is a surjective and covering map (see [WG]). Since, moreover,
$S^3-S$ is simply connected (by Proposition~\ref{transv}),
$G:M-G^{-1}(S)\rightarrow S^3-S$ is a diffeomorphism. To prove that $M^3$ is
homeomorphic to $S^3$, it is necessary and sufficient that every piecewise
smooth simple closed
curve $\ga\subset M^3$ is contained in a topological cube $Q\subset M^3$ (by
Lemma~\ref{Bing}).

In order to prove that every piecewise smooth curve $\ga$ lies in a topological cube,
observe that we may suppose that $\ga \cap K=\emptyset$ (here $K=G^{-1}(S)$). Indeed, by lemma
\ref{difeo} there exists a diffeomorphism $h$ of $M$ such $h(\ga)\cap K= \emptyset.$ Then, if
$h(\ga)$ lies in a topological cube $R$, the $\ga$ itself lies in the
topological cube $h^{-1}(R)$ too, thus we can, in fact, make this assumption.

 Now, since $\ga\subset M-K$ and $M-K$ is diffeomorphic to $S^3-S$, we
may consider $\ga\subset \real^3-S$, $S$ a compact subset of $\real^3$ with
Hausdorff dimension less than $1$ via identification by the diffeomorphism $G$
and stereographic projection. In this case, we can follow the proof of
Proposition~\ref{transv} to conclude that $f:\ga\times S\rightarrow S^2$,
$f(x,y)=(x-y)/||x-y||$ is Lipschitz. Because $\overline{\dim_B\ga}\leq 1, \dim_H
S<1$ (here we are using that $\ga$ is piecewise smooth), we obtain a direction $v\in S^2$ such that
$F := \bigcup\limits_{t\in\real} (L_t(\ga))$ is disjoint from $S$, where $L_t(p):=p+t\cdot
v$. By compactness of $\ga$ it is easy that $F$ is a closed subset of $\real^n$. This implies
that $3\ \ep=d(F,S)>0$. Consider $F_{\ep}=\{x:d(x,F)\leq\ep\}$ and $S_{\ep}=\{x:
d(x,S)\leq\ep\}$. By definition of $\ep>0$, $F_{\ep}\cap S_{\ep}=\emptyset$,
then we can choose $\vfi:\real^3\rightarrow \real$ a smooth function such that
$\vfi|_{F_{\ep}} = 1$, $\vfi|_{S_{\ep}} = 0$. Consider the vector field $X(p)=\vfi(p)\cdot v$
and let $X_t$, $t\in\real$ the $X$-flow. We have $X_t(p)
=p+tv\ \forall p\in\ga$ and $X_t(p)=p\ \forall p\in S$, for any $t\in\real$.
Choosing $N$ real such that $S\subset B_N(0)$ and $T$ such that $t\geq T\Rightarrow
L_t(\ga)\cap B_N(0)=\emptyset$, we obtain a global homeomorphism $X_t$ which sends $\ga$
outside $B_N(0)$ and keep fixed $S$, $\forall \ t\geq T$.

Observe that $X_t(\ga)$ is
contained in the interior of a topological cube $Q\subset \real^3-B_N(0)$. Then,  observing that $X_t$ is a
diffeomorphism  and that $X_t(x)=x$ for every $x\in S$ and $t\in \real$, we have that $\ga \subset
X_{-t}(Q) \subset \real^3-S$, $\forall\ t\geq T$. This concludes the proof.
\end{proof}

\section{Proof of Theorem A in the case $n\geq 4$}
We start this section with the statement of generalized Poincar\'e Conjecture :

\begin{theorem}\label{Poinc}A compact simply connected homological sphere $M^n$ is
homeomorphic to $S^n$, if $n\geq 4$ (diffeomorphic for $n=5,6$).
\end{theorem}

The proof of generalized Poincar\'e Conjecture is due to Smale~\cite{Smale} for
$n\geq 5$ and to Freedman~\cite{Freedman} for $n=4$. This theorem makes the proof of
the Theorem B a little bit easier than the proof of Theorem A.

\begin{proof}[Proof of Theorem A in the case $n\geq 4$]If $k=n$, there is nothing to prove.
Indeed, in this case, $G:M^n\rightarrow S^n$ is a diffeomorphism, by definition.
I.e., without loss of generality we can suppose $k\leq n-1$; $\al_k(M)<
k-[\frac{n}{2}]\Rightarrow \exists\ x:M^n\rightarrow \real^{n+1}$ immersion,
$rank(x)\geq k$, $dim_H(x)<k-[\frac{n}{2}]$. The hypothesis implies that $M-G^{-1}(S)$ is
connected,
$S^n-S$ is simply connected  and $G$ is a proper map whose jacobian never
vanishes. By~\cite{WG}, $G$ is a surjective, covering map. So, we conclude that $G:M-G^{-1}(S)\rightarrow
S^n-S$ is diffeomorphism. But $S^n-S$ is $(n-1-k+[\frac{n}{2}])$-connected, by
Proposition~\ref{transv}. In particular, because $k\leq n-1$, $S^n-S$ is
$[\frac{n}{2}]$-connected and so, using the diffeomorphism $G$, $M-K$ is
$[\frac{n}{2}]$-connected, where $K=G^{-1}(S)$. It is sufficient to prove that
$M^n$ is a simply connected homological sphere, by Theorem~\ref{Poinc}. By Lemma~\ref{difeo}, $M-K$ is $[\frac{n}{2}]$-connected and
$dim_H(K)<n-[\frac{n}{2}]$ (by Proposition\ref{transv}) implies $M$ itself is $[\frac{n}{2}]$-connected.
It is know that $H^i(M)=L(H_i(M))\oplus T(H_{i-1}(M))$, L and T denotes the free
part and the torsion part of the group. By Poincar\'e duality, $H_{n-i}(M)\simeq
H^i(M)$. The fact that $M$ is $[\frac{n}{2}]$-connected and the other facts give us
$H_i(M)=0, \text{ for } 0<i<n$. This concludes the proof.
\end{proof}

\section{Proof of Theorems B and C}
In this section we make some comments on extensions of Theorem~\ref{theoremA}.
Although these extensions are quite easy, they
were omitted so far to make the presentation of the paper more clear. Now, we
are going to improve our previous results. First, all
preceding arguments works with assumption that $H_{k-[\frac{n}{2}]}(S)=0$ and
$rank(x)\geq k$ in Theorems A, B (where $H_s$ is the $s$-dimensional Hausdorff
measure). We prefer to consider the hypothesis as its
stands in these theorems because it is more interesting to define the invariants $\al_{k}(M)$. The
reason to this ``new'' hypothesis works is that our proofs, essentially, depend on the
existence of special directions $v\in S^{n-1}$. But these directions exist if
the singular sets have Hausdorff measure 0. Secondly, $M$ need not to be compact.
It is sufficient that $M$ is of \emph{finite geometric type} (here our definition of finite
geometrical type is a little bit different from \cite{BFM}). We will make more
precise these comments in proof of Theorem 6.2 below, after recalling the definition :

\begin{definition}Let $\overline{M}^n$ a compact (oriented) manifold and
$q_1,\dots,q_k\in\overline{M}^n$. Let $M^n=\overline{M}^n-\{q_1,\dots,q_k\}$. An
immersion $x:M^n\rightarrow \real^{n+1}$ is of finite geometrical type if $M^n$
is complete in the induced metric, the Gauss map $G:M^n\rightarrow S^n$ extends
continuously to a function $\overline{G}:\overline{M}^n\rightarrow S^n$ and the set
$G^{-1}(S)$ has $H_{n-1}(G^{-1}(S))=0$ (this last condition occurs
if $rank(x)\geq k$ and $\dim_H(x)<k-1$, or more generally, if $rank(x)\geq k$
and $H_{k-1}(S)=0$).
\end{definition}

As pointed out in the introduction, the conditions in the previous definition are satisfied, for example, by complete
hypersurfaces with finite total curvature whose Gauss-Kronecker curvature $H_n=
k_1\dots k_n$ does not change of sign and vanish in a small set, as showed by \cite{dCE}. Recall that
a hypersurface $x:M^n\rightarrow \real^{n+1}$ has total finite curvature if
$\int_M |A|^n dM<\infty$, $|A|=(\sum\limits_i k_i^2)^{1/2}$, $k_i$ are the
principal curvatures. Then, there are examples satisfying the
definition. With these observations, we now prove our Theorem~\ref{theoremB}

\begin{proof}[Proof of Theorem~\ref{theoremB}]To avoid unnecessary repetitions, we will only indicate the
principal modifications needed in proof of Theorems~\ref{theoremA} and~\ref{theoremB} by stating ``new''
propositions, which are analogous to the previous ones, and making a few comments
in their proofs. The details are left to reader.

\begin{proposition}[proposition~\ref{dim pre imagem}']\label{dim pre imagem'}Let
$f:M^m\rightarrow N^n$ a $C^1$-map and $A\subset N$. If $H_s(A)=0$, then
$H_{s+n-rank(f)}(f^{-1}(A))=0$.
\end{proposition}

\begin{proof}It suffices to show that for any $p\in f^{-1}(A)$, there is an
open set $U=U(p)\ni p$ such that $H_{s+n-r}(f^{-1}(A)\cap U)=0$. However, if $U$ is
chosen as in proof of Proposition~\ref{dim pre imagem}, we have $f^{-1}(A)\cap
U\subset \vfi(\pi(A)\times\real^{n-r})$, where $\vfi$ is a diffeomorphism, $r=
rank(f)$ and $\pi$ is the projection in first $r$ variables. By
Propositions~\ref{fatos basicos},~\ref{produto}, $H_{s+n-r}(f^{-1}(A)\cap U)\leq
C_1\cdot
H_{s+n-r}(\pi(A)\times\real^{n-r})\leq C_1\cdot C_2\cdot H_s(A) = 0$, where
$C_1$ depends only on $\vfi$ and $C_2$ depends only on $(n-r)$. This finishes
the proof.
\end{proof}

\begin{proposition}[Proposition~\ref{transv}']\label{transv'}Let $n\geq 3$ and $F$ a
closed subset of $M^n$ such that $H_{n-1}(F)=0$ then $M-F$ is connected. If $M^n=
\real^n$ or $M^n=S^n$, $F$ is compact and $H_{n-k}(F)=0$ then $M^n-F$ is
$k$-connected.
\end{proposition}

\begin{proof}First, if $\ga$ is a path in $M^n$ from $x$ to $y$, $x,y\notin F$,
we can suppose $\ga$ a smooth simple curve. So, there is a compact tubular
neighborhood $\SL=\ga\times D^{n-1}$ of $\ga$. Since $dim(\SL_p)=n-1$, $F\cap\SL_p$ has Lebesgue
measure 0 for any $p$. Thus, using that $F$ is closed and $\ga$ is compact, we
obtain some arbitrarily small vector $v$ such that $(\ga\times v)\cap F=\emptyset$.
Then, $M^n-F$ is connected.

Second, if $[\Ga]\in\pi_i(\real^n-F),\ i\leq k$ is a homotopy class and $\Ga$ is
a smooth representative, define $f:\Ga\times F\rightarrow S^{n-1}$,
$f(x,y)=(x-y)/||x-y||$. Following the proof of Proposition~\ref{transv}, $f$ is
Lipschitz. Now, since $H_{n-k-1}(F)=0$, we have, by Proposition~\ref{produto},
$H_{n-1}(\Ga\times F)=0$. Thus, Proposition~\ref{fatos basicos} implies
$H_{n-1}(f(\Ga\times F))=0$. This concludes the proof.
\end{proof}

\begin{lemma}[Lemma~\ref{difeo}']\label{difeo'}Suppose that $\Ga\in\pi_i(M^n)$
is Lipschitz and $K\subset M^n$ is compact, $H_{n-i}(K)=0$. Then there are
diffeomorphisms $h$ of $M$, arbitrarily close to identity map, such that $h(\Ga) \cap
K=\emptyset$.
\end{lemma}

\begin{proof}If $\Ga$ is Lipschitz and $\Ga$ lies in a parametrized
neighborhood, we can take $F:\Ga\times K\rightarrow\real^n$, $F(x,y)=x-y$ a
Lipschitz function. Because $H_n(\Ga\times K)=0$, then $H_n(F(\Ga\times
K))=0$. In the general case we proceed as in proof of Lemma~\ref{difeo}. Take, by
compactness, a finite number of parametrized neighborhoods and apply the previous
case. By finiteness of number of parametrized neighborhoods and using that $K$ is compact, an induction argument achieves
the desired diffeomorphisms $h$. This
concludes the proof.
\end{proof}

Returning to proof of Theorem~\ref{theoremB}, observe that in Theorem~\ref{theoremA}, we
need $\overline{G}:\overline{M}^n-\overline{G}^{-1}(\widetilde{S})\rightarrow S^n-\widetilde{S}$ to be a
diffeomorphism, where $\widetilde{S}=S\cup\{\overline{G}(q_i): i=1,\dots,k\}$. This remains true
because $(*) \ H_{k-[\frac{n}{2}]}(S)=0$ implies $S^n-\widetilde{S}$ is
$(n-1-k+[\frac{n}{2}])$-connected. In fact, this is a consequence of $(*)$,
Proposition~\ref{transv'} and $\{p_i : i=1,\dots,k\}$ is finite ($p_i:= \overline{G}(q_i)$). Moreover, $rank(x)\geq k$ imply, by Proposition~\ref{dim pre
imagem'},~\ref{transv'}, $\overline{M}-G^{-1}(\widetilde{S})$ is connected. Indeed, these
propositions says that $rank(x)\geq k \Rightarrow H_{n-[\frac{n}{2}]}(G^{-1}(S))=0$ and
$H_{n-1}(G^{-1}(S))=0\Rightarrow M-G^{-1}(S)\text{ is connected}$. However, if
$\overline{G}^{-1}(\widetilde{S})-( G^{-1}(S)\cup\{q_i : i=1,\dots,k \}):= A $, then, for all
$x\in A$, $(**)\ \det D_xG\neq 0$. In particular, since $G(A)\subset\{p_i :
i=1,\dots,k\}$,
$(**)$ implies $\dim_H(A)=0$. Then, $H_{n-[\frac{n}{2}]}(\overline{G}^{-1}(\widetilde{S})) =
H_{n-[\frac{n}{2}]}(G^{-1}(S))=0$. Thus, by~\cite{WG}, G is
surjective and covering map (because it is proper and its jacobian never
vanishes). In particular, by simple connectivity, $G$ is a diffeomorphism. At
this point, using the previous lemma and propositions, it is sufficient to follow
proof of Theorem~\ref{theoremA}, if $n=3$, and proof of Theorem~\ref{theoremB}, if $n\geq 4$, to obtain
$\overline{M}^n\simeq S^n$. This concludes the proof.
\end{proof}

For even dimensions, we can follow \cite{BFM} and improve Theorem~\ref{theoremB} :

\begin{theorem}[Theorem~\ref{theoremC}]Let $x:M^{2n}\rightarrow\real^{2n+1}$, $n\geq 2$ an immersion of finite
geometrical type with $H_{k-n}(S)=0$, $rank(x)\geq k$. Then $M^{2n}$ is
topologically a sphere minus two points. If $M^{2n}$ is minimal, $M^{2n}$ is a
$2n$-catenoid.
\end{theorem}

For sake of completeness we present an outline of proof of Theorem~\ref{theoremC}.

\begin{proof}[Outline of proof of Theorem~\ref{theoremC}]Barbosa, Fukuoka, Mercuri define to each \emph{end} $p$ of
$M$ a geometric index $I(p)$ that is related with the topology of $M$ by the
formula (see theorem 2.3 of \cite{BFM}):
\begin{equation}\chi(\overline{M}^{2n})=\sum\limits_{i=1}^{k}(1+I(p_i))+2\si m
\end{equation}
where $\si$ is the sign of Gauss-Kronecker curvature and $m$ is the degree of
$G:M^n\rightarrow S^n$. Now, the hypothesis $2n>2$ implies (see~\cite{BFM})
$I(p_i)=1,\forall\ i$. Since we know, by theorem 6.2, $\overline{M}^{2n}$ is a
sphere, we have $2=2k+2\si m$. But, it is easy to see that $m=deg(G)=1$ because $G$ is
a diffeomorphism outside the singular set. Then, $2=2k+2\si\Rightarrow k=2, \
\si=-1$. In particular, $M$ is a sphere minus two points.

If $M$ is minimal, we will use the following theorem of Schoen : \textit{The
only minimal immersions, which are regular at infinity and have two ends, are
the catenoid and a pair of planes}. The regularity at infinity in our case holds
if the ends are embedded. However, $I(p)=1$ means exactly this. So, we can use
this theorem in the case of minimal hypersurfaces of finite geometric type. This
concludes the
outline of proof.
\end{proof}

\begin{remark}We can extend theorem A in a different direction (without
mention of $rank(x)$). In fact, using only that $G$ is Lipschitz, it suffices assume that
$H_{n-[\frac{n}{2}]}(C)=0$ ($C$ is the set of points where the Gauss-Kronecker
curvature vanishes). This is essentially the hypothesis of Barbosa, Fukuoka and
Mercuri~\cite{BFM}. We prefer to state Theorems~\ref{theoremB} and~\ref{theoremC} as before since the
classical theorems concerning estimates for Hausdorff dimension (Morse-Sard,
Moreira) deal only with the critical values $S$ and, in particular, our
Corollary 2.4 will be more difficult if the hypothesis is changed to $H_1(C)=0$
for some immersion $x:M^3\rightarrow\real^4$ (although, in this assumption, we
have no problems with $rank(x)$, i.e., this assumption has some advantages).
\end{remark}

\begin{remark}It would be interesting to know if there are examples of codimension $1$
immersions with singular set which is not in the situation of
Barbosa-Fukuoka-Mercuri and do Carmo-Elbert but instead satisfies our hypothesis.
This question was posed to the second author by Walcy Santos during the
Differential Geometry
seminar at IMPA.
In fact, these immersions can be constructed with some extra work. Some examples will be
presented in another work to appear elsewhere.
\end{remark}

\section{Final Remarks}\label{finalremarks}

The  Corollary~\ref{conjecture1} is motivated by Anderson's program for Poincar\'e
Conjecture. In order to coherently describe this program, we briefly recall some facts
about topology of $3$-manifolds.

An attempt to better understand the topology of $3$-manifolds (in particular,
give an answer to Poincar\'e Conjecture) is the so called ``Thurston
Geometrization Conjecture''. Thurston's Conjecture goes beyond Poincar\'e
Conjecture (which is a very simple corollary of this conjecture). In fact, its goal is the understanding of 3-manifolds by decomposing
them into pieces which could be ``geometrizated'', i.e., one could put complete
locally homogeneous metric in each of this pieces. Thurston showed that, in
three dimensions, there are exactly \emph{eight} geometries, all of which are
realizable. Namely, they are : the constant curvature spaces $\mathbb{H}^3$,
$\real^3$, $S^3$, the products $\mathbb{H}^2\times \real$, $S^2\times\real$ and
the twisted products $\widetilde{SL}(2,\real)$, Nil, Sol (for details see \cite{Thurston}). Thurston proved his conjecture in some particular cases (e.g.,
for \emph{Haken} manifolds). These particular cases are not easy. To prove the
result Thurston developed a wealth of new geometrical ideas and machinery to carry
this out. In a few words, Thurston's proof is made by induction. He decomposes the
manifold $M$ in an appropriate hierarchy of submanifolds $M_k=M\supset\dots\supset\text{union of balls}=
M_0$ (this is possible if $M$ is
\emph{Haken}). Then, if $M_{i-1}$ has a metric with some
properties, it is possible glue certain ends of $M_{i-1}$ to obtain $M_i$.
Moreover, by a deformation and isometric gluing of ends argument, $M_i$ has a metric with the same properties
of that from $M_{i-1}$. This is the most difficult part of the proof. So, the
induction holds and $M$ itself satisfies the Geometrization Conjecture.

Recently, M.~Anderson \cite{Anderson} formulated three conjectures that imply
Thurston's Conjecture. Morally, these three
conjectures says that information about the \emph{sigma} constant give us
information about geometry and topology of 3-manifolds. We recall the definition
of sigma constant. If $S(g)
:=\int_M s_g\ dV_g$ is the total scalar curvature functional ($g$ is a metric
with unit volume, i.e., $g\in\SM_1$, $dV_g$ is volume form determined by $g$ and $s_g$ is the scalar
curvature) and $[g]:=\{\widetilde{g}\in\SM_1:\widetilde{g}=\psi^2 g \text{, for
some smooth positive function $\psi$}\}$ is the conformal class of $g$, then $S$ is a bounded from below functional in
$[g]$. Thus, we can define $\mu[g]=\inf\limits_{g\in[g]} S(g)$ called
\emph{Yamabe} constant of $[g]$. An elementary comparison argument shows
$\mu[g]\leq \mu(S^n,g_{can})$, where $g_{can}$ is the canonical metric of $S^n$
with unit 1 and positive constant curvature. Then makes sense to define the
sigma constant :
\begin{equation}\label{constante sigma}\si(M)=\sup\limits_{[g]\in\SC}\mu[g]
\end{equation} where $\SC$ is the space of all conformal classes. The sigma
constant is a smooth invariant defined by a minimax principle (see
equation~\ref{constante sigma}). The first part of this minimax procedure was solved by
Yamabe~\cite{Yamabe}. More precisely, for any conformal class $[g]\in\SC$,
$\mu[g]$ is realized by a (smooth) metric $g_{\mu}\in[g]$ such that
$s_{g_{\mu}}\equiv \mu[g]$ (a such $g_{\mu}$ is called \emph{Yamabe} metric).
The second part of this procedure is more difficult since it depends on the
underlying topology. The sigma constant is important since is know that critical points of the
scalar curvature functional $S$ are Einstein metrics. But is not know if
$\si(M)$ is a critical value of $S$ (partially by non-uniqueness of Yamabe
metrics). Then, if one show that is possible to realize the second part of
minimax procedure and that $\si(M)$ is a critical value of $S$, we obtain the
Geometrization Conjecture.

The approach above is very difficult. To see this, we remark that all of three Anderson's Conjectures are
necessary to obtain the ``Elliptization Conjecture'' (the particular case of
Thurston's Conjecture which implies Poincar\'e Conjecture). In others words, we
have to deal with all cases of Thurston Conjecture to obtain Poincar\'e
Conjecture. This inspires our definition of other minimax smooth invariants.
The advantage in these invariants is that they do not
requires construction of metrics with positive constant curvature. But the
disadvantage is we always
work extrinsically.

To finish the paper, we comment that there are many others attacks and
approachs to Poincar\'e Conjecture. For example, see~\cite{Gabai} for an
accessible exposition of V.~Po\'enaru's program and~\cite{Poenaru} for recent
proof of one step of this program. In the other hand, some authors (e.g., Bing
~\cite{Bing}) believes that only simple connectivity is not sufficient for a
manifold be $S^3$. Recentely, Perelman  has proposed a proof of the geometrization conjecture.

\textbf{Acknowledgments.} The authors are thankful  to Jairo Bochi, Jo\~ao
Pedro Santos, Carlos Morales, Alexander Arbieto and Carlos Moreira for fruitful conversations.  Also,
the authors are indebted to Professor Manfredo do Carmo for his kind encouragement and to
professor Marcelo Viana for his suggestions and advice. The comments of
anonymous referees were useful to improve the presentation of this paper and correct several misspells. Last, but not least, we are thankful to IMPA and its staff.


\vspace{1cm}

\noindent Carlos Matheus ( matheus{\@@}impa.br )\\
          IMPA, Estrada Dona castorina, 110 Jardim Bot\^anico 22460-320\\
          Rio de Janeiro, RJ, Brazil\\
          Krerley Oliveira ( krerley{\@@}mat.ufal.br )\\
UFAL, Campus AC Simoes,s/n Tabuleiro, 57072-090 \\
Maceio, AL, Brazil

\end{document}